\documentclass[12pt,a4paper]{amsart}

\newtheorem{Result}{Result}[section]
\newtheorem{Thm}  [Result]{Theorem}
\newtheorem{Prop} [Result]{Proposition}

\newtheorem{Cor}  [Result]{Corollary}

\newtheorem*{Def}{Definition}
\newtheorem*{Rem}{Remark}

\def\qedbox{\hbox{\vrule\vbox{\hrule width 6pt\vskip6pt\hrule}\vrule}}
\def\qed{\ifvmode\leavevmode\fi
  \unskip\nobreak\hfill\penalty50 \quad \null\nobreak\hfill
  \qedbox{\parfillskip0pt \finalhyphendemerits0 \par}}
\newenvironment{Proof}{\noindent{\it Proof.~}\ignorespaces}{\qed\medskip}

\def\rnbox#1{(\romannumeral#1)}
\def\itemn#1{\item[\rnbox{#1}]}

\def\N{{\mathbb{N}}}

\def\ptensor{\mathop{\hat\otimes_\pi}}
\def\itensor{\mathop{\hat\otimes_\varepsilon}}

\def\JH{\hbox{\it JH\/}}

\def\map#1#2#3{#1\colon#2\,{\longrightarrow}\,#3}

\begin{document}

\title[Subprojectivity and superprojectivity of Banach spaces]
  {On subprojectivity and superprojectivity of Banach spaces}

\author[E.M. Galego, M. Gonz\'alez and J. Pello]{El\'oi M. Galego,
    Manuel Gonz\'alez and Javier Pello}

\address{Departamento de Matem\'atica,
    Instituto de Matem\'atica e Estat\'\i stica,
    Universidade de S\~ao Paulo, 
    \hbox{05508-090} S\~ao Paulo SP, Brazil.}
\email{eloi@ime.usp.br}

\address{Departamento de Matem\'aticas, Facultad de Ciencias,
    Universidad de Cantabria, \hbox{E-39071} Santander, Spain}
\email{manuel.gonzalez@unican.es}

\address{
    Escuela Superior de Ciencias Experimentales y Tecnolog\'\i a,
    Universidad Rey Juan Carlos, \hbox{E-28933} M\'ostoles, Spain}
\email{javier.pello@urjc.es}

\thanks{Supported in part by MICINN (Spain), Grant MTM2013-45643.\\
    2010 Mathematics Subject Classification. Primary: 46B20, 46B03.\\
    Keywords: Banach space; subprojective space; superprojective space.}

\begin{abstract}
We obtain some results for and further examples of subprojective and
superprojective Banach spaces.
We also give several conditions providing examples of non-reflexive
superprojective spaces; one of these conditions is stable under $c_0$-sums
and projective tensor products.
\end{abstract}

\maketitle

\begin{section}{Introduction}

The classes of subprojective and superprojective Banach spaces were
introduced by Whitley \cite{whitley} to find conditions for the conjugate of
an operator to be strictly singular or strictly cosingular.
They are relevant in the study of the perturbation classes problem for
semi-Fredholm operators \cite{gonzalez}, which has a positive answer when
one of the spaces is subprojective or superprojective \cite{gonzalez-et-al}.
A reflexive space is subprojective (superprojective) if and only if
its dual is superprojective (subprojective).
In general, however, $X$ being subprojective does not imply that $X^*$ is
superprojective, and $X^*$ being subprojective does not imply that $X$ is
superprojective, and it is unknown whether the remaining implications
are valid \cite[Introduction]{gonzalez-pello}.
Basic examples of subprojective spaces are $\ell_p$ for $1\leq p < \infty$
and $L_p(0,1)$ for $2\leq p<\infty$ \cite[Proposition 2.4]{gonzalez-et-al};
and $C(K)$ spaces with $K$ a scattered compact are both subprojective and
superprojective \cite[Propositions 2.4 and 3.4]{gonzalez-et-al}. Moreover,
recent systematic studies of subprojective spaces \cite{oikhberg-spinu}
(see also \cite{galego-samuel}) and superprojective spaces
\cite{gonzalez-pello}
have widely increased the family of known examples in those classes.

Here we continue the study of subprojective and superprojective Banach
spaces. In Section~\ref{sect-char} we give some characterisations of these
classes of spaces in terms of improjective operators, and apply them to
analyse the subprojectivity and superprojectivity of spaces with the
Dunford-Pettis property, in particular $\mathcal{L}_1$-spaces and
$\mathcal{L}_\infty$-spaces. We show that hereditarily-$\ell_1$ spaces
with an unconditional basis and hereditarily-$c_0$ spaces are subprojective,
and that $C([0,\lambda],X)$ is subprojective when $X$ is subprojective and
$\lambda$ is an arbitrary ordinal. We also study the subprojectivity and
superprojectivity of some $\mathcal{L}_\infty$-spaces obtained by Bourgain
and Delbaen~\cite{bourgain-delbaen}, which provide counterexamples to some
natural conjectures.

In Section~\ref{sect-cond} we find new examples of non-reflexive
superprojective Banach spaces.
We show that, if $X$ has property~(V) and $X^*$ is hereditarily~$\ell_1$,
then $X$ is superprojective.
In particular, this is true for the spaces in the class that we denote
by $Sp(U^{-1}\circ W)$, which includes~$C(K)$ spaces with $K$ scattered, the
isometric preduals of~$\ell_1(\Gamma)$ and Hagler's space~$\JH$~\cite{hagler}.
Note that $\JH$ is a separable space that contains no copies of~$\ell_1$ and
has non-separable dual, hence $\JH$ does not admit an unconditional basis.
The class $Sp(U^{-1}\circ W)$ is shown to be stable under passing to
quotients and under taking projective tensor products and $c_0$-sums.
We also show that the predual $d(w,1)_*$ of the Lorentz space $d(w,1)$ and
the Schreier space~$S$ are superprojective, although they do not belong to
$Sp(U^{-1}\circ W)$, and that their dual spaces are subprojective,
but the tensor products $S\ptensor S$ and $S\ptensor \ell_p$ 
are not superprojective.

\medskip

In the sequel, subspaces of a Banach space are assumed to be closed unless
otherwise stated.
Given a subspace~$M$ of a Banach space, $J_M$ and~$Q_M$ denote its
natural embedding and quotient map.
A Banach space~$X$ is hereditarily~$Z$ if every infinite-dimensional
subspace of~$X$ contains a subspace isomorphic to~$Z$.
Given Banach spaces $X$ and~$Y$, $L(X,Y)$ denotes the set of all (continuous,
linear) operators from $X$ into~$Y$, and $K(X,Y)$ denotes the subset of
compact operators.

An injection is an isomorphic embedding with infinite-dimensional range,
and a surjection is a surjective operator with infinite-dimensional range.
A compact space~$K$ is said to be scattered, or dispersed, if every nonempty
subset of~$K$ has an isolated point.

A Banach space $X$ is an $\mathcal{L}_{p,\lambda}$-space
($1\leq p\leq\infty$; $1\leq \lambda <\infty$) if every finite-dimensional
subspace~$F$ of~$X$ is contained in another finite-dimensional subspace~$E$
of~$X$ whose Banach-Mazur distance to the space $\ell^{\dim E}_p$ is
at most~$\lambda$.
The space~$X$ is an $\mathcal{L}_p$-space if it is an
$\mathcal L_{p,\lambda}$-space for some~$\lambda$.

\end{section}

\begin{section}{Subprojective and superprojective spaces}
\label{sect-char}

We begin by recalling the definitions given in \cite{whitley} of the concepts
we investigate.

\begin{Def}
A Banach space~$X$ is called \emph{subprojective} if every
infinite-dimensional subspace of~$X$ contains an infinite-dimensional
subspace complemented in~$X$, and $X$ is called \emph{superprojective}
if every infinite-codimensional subspace of~$X$ is contained in an
infinite-codimensional subspace complemented in~$X$.
\end{Def}

The following result \cite[Proposition 3.3]{gonzalez-pello} is useful
to show that some spaces fail subprojectivity or superprojectivity.

\begin{Prop}\label{ell-1}
If a Banach space~$X$ contains a copy of~$\ell_1$, then $X$ is not
superprojective and $X^*$ is not subprojective.
\end{Prop}

An operator $\map{T}{X}{Y}$ is called \emph{strictly singular} if there is
no infinite-dimensional subspace~$M$ of~$X$ such that the restriction $TJ_M$
is an isomorphism. The following, more general concept was introduced by
Tarafdar \cite{tarafdar}.

\begin{Def}
An operator $\map{T}{X}{Y}$ is called \emph{improjective} if there is no
infinite-dimensional subspace~$M$ of~$X$ such that the restriction $TJ_M$ is
an isomorphism and $T(M)$ is complemented in~$Y$.
\end{Def}

An operator $\map{T}{X}{Y}$ is called \emph{strictly cosingular} if there is
no infinite-codimensional subspace~$N$ of~$Y$ such that $Q_NT$ is surjective.
The following characterisation, obtained in
\cite[Theorem 2.3]{aiena-gonzalez},
shows that strictly cosingular operators are improjective.

\begin{Prop}\label{improj}
An operator $\map{T}{X}{Y}$ is improjective if and only if there is no
infinite-codimensional subspace~$N$ of~$Y$ such that $Q_NT$ is surjective
and $T^{-1}(N)$ is complemented in~$X$.
\end{Prop}

Next we give some characterisations of subprojectivity and superprojectivity
in terms of improjective operators.

\begin{Prop}\label{subproj-char}
For a Banach space~$X$ the following are equivalent:
\begin{enumerate}
\itemn1 $X$ is subprojective;
\itemn2 every improjective operator $\map{T}{Z}{X}$ is strictly singular;
\itemn3 there exists no improjective injection $\map{J}{Z}{X}$.
\end{enumerate}
\end{Prop}

\begin{Proof}
\rnbox1 $\Rightarrow$ \rnbox2
Suppose that $X$ is subprojective and an operator $\map{T}{Z}{X}$
is not strictly singular.
Then there exists an infinite-dimensional subspace~$M$ of~$Z$ such that
$TJ_M$ is an isomorphism.
Let $N$ be an infinite-dimensional subspace of $T(M)$ complemented in~$X$;
then $T$ is an isomorphism on $M_0 := M\cap T^{-1}(N)$ and $T(M_0) = N$,
hence $T$ is not improjective.

\rnbox2 $\Rightarrow$ \rnbox3 It is enough to observe that injections
are not strictly singular.

\rnbox3 $\Rightarrow$ \rnbox1 Given an infinite-dimensional subspace~$M$
of~$X$, the injection $\map{J_M}{M}{X}$ is not improjective, so there exists
an infinite-dimensional subspace~$N$ of~$M$ which is complemented in~$X$.
Thus $X$ is subprojective.
\end{Proof}

\begin{Prop}\label{superproj-char}
For a Banach space $X$ the following are equivalent:
\begin{enumerate}
\itemn1 $X$ is superprojective;
\itemn2 every improjective operator $\map{T}{X}{Y}$ is strictly cosingular;
\itemn3 there exists no improjective surjection $\map{Q}{X}{Y}$.
\end{enumerate}
\end{Prop}

\begin{Proof}
\rnbox1 $\Rightarrow$ \rnbox2
Suppose that $X$ is superprojective and an operator $\map{T}{X}{Y}$
is not strictly cosingular.
Then there exists an infinite-codimensional subspace~$N$ of~$Y$ such that
$Q_NT$ is surjective.
Let $M$ be an infinite-codimensional subspace complemented in~$X$ and
containing $T^{-1}(N)$; then $T(M)$ is closed and infinite-codimensional,
$Q_{T(M)}T$ is surjective and $T^{-1}T(M) = M$ is complemented in~$X$,
hence $T$ is not improjective by Proposition~\ref{improj}.

\rnbox2 $\Rightarrow$ \rnbox3 It is enough to observe that surjections
are not strictly cosingular.

\rnbox3 $\Rightarrow$ \rnbox1 Given an infinite-codimensional subspace~$N$
of~$X$, the surjection $\map{Q_N}{X}{X/N}$ is not improjective, so there
exists an infinite-codimensional subspace~$M$ containing~$N$ which is
complemented in~$X$. Thus $X$ is superprojective.
\end{Proof}

A Banach space $X$ has the Dunford-Pettis property (DPP in short) if every
weakly compact operator $\map{T}{X}{Y}$ takes weakly convergent sequences
to convergent sequences; or, equivalently, if every weakly compact operator
$\map{T}{X}{Y}$ takes weakly compact sets to relatively compact sets.
We refer the reader to \cite[Section~5.4]{kalton-albiac} and
\cite[Section~10]{handbook} for further information on the DPP.
Examples of spaces with the DPP are the $\mathcal{L}_\infty$-spaces
and the $\mathcal{L}_1$-spaces \cite[Section~10]{handbook}; in particular,
the spaces of continuous functions on a compact $C(K)$ and the spaces of
integrable functions $L_1(\mu)$.

The next result establishes some necessary conditions for spaces with the DPP
to be subprojective or superprojective.

\begin{Prop}\label{negative-DPP}
Let $X$ be a Banach space satisfying the DPP.
\begin{enumerate}
\item If $X$ is subprojective, then it contains no infinite-dimensional
reflexive subspaces.
\item If $X$ is superprojective, then it admits no infinite-dimensional
reflexive quotients.
\end{enumerate}
\end{Prop}

\begin{Proof}
(1) Let $R$ be a reflexive subspace of~$X$. By Proposition~\ref{subproj-char},
it is enough to show that the embedding $\map{J_R}{R}{X}$ is strictly
cosingular, hence improjective, as that would make $R$ finite-dimensional.

Let $\map{Q}{X}{Z}$ be an operator such that $QJ_R$ is surjective.
Then $QJ_R$ is weakly compact, so $Z$ is reflexive and $Q$ itself is weakly
compact, hence completely continuous by the DPP of $X$.
Thus $QJ_R$ is compact, and $Z$ is finite-dimensional.

(2) We could apply Proposition~\ref{superproj-char} to give a proof similar
to that of (1), but we choose an alternative one.
Take a bounded sequence $(x_n)_{n\in\N}$ in~$X$ whose image in the reflexive
quotient~$X/M$ is weakly convergent but does not have any convergent
subsequences. Then $Q_M$ is weakly compact and $X$ has the DPP, so
$Q_M$ takes weakly Cauchy sequences to convergent sequences and
$(x_n)_{n\in\N}$ cannot have any weakly Cauchy subsequence.
Thus $X$ contains a subspace isomorphic to~$\ell_1$ and it is not
superprojective by Proposition~\ref{ell-1}.
\end{Proof}

\begin{Cor}\label{negative-DPP-L1}
A $\mathcal{L}_1$-space is subprojective if and only if it contains no
infinite-dimensional reflexive subspaces.
\end{Cor}

\begin{Proof}
The direct implication is a consequence of Proposition~\ref{negative-DPP}.
For the converse, observe that each $\mathcal{L}_1$-space $X$ is isomorphic
to a subspace of some space $L_1(\mu)$ \cite{lindenstrauss-pelczynski}.
Therefore, every non-reflexive subspace of~$X$ contains a copy of~$\ell_1$
complemented in~$X$ \cite[Proposition~5.6.2]{kalton-albiac}.
\end{Proof}

The analogue of Corollary~\ref{negative-DPP-L1} for
$\mathcal{L}_\infty$-spaces does not hold.
We will see later that there exists a $\mathcal{L}_\infty$-space $Y_{bd}$
admitting no infinite-dimensional reflexive quotient which is not
superprojective.

\medskip

The next result was essentially proved by D\'\i az and
Fern\'andez~\cite{diaz-fernandez}.

\begin{Prop}\label{hered-c0}
Every hereditarily-$c_0$ Banach space is subprojective.
\end{Prop}

\begin{Proof}
It was proved in \cite[Theorem 2.2]{diaz-fernandez} that if a Banach
space~$X$ contains no copies of~$\ell_1$, then every copy of~$c_0$ in~$X$
contains another copy of~$c_0$ which is complemented in~$X$.
\end{Proof}

There are hereditarily-$c_0$ spaces that admit~$\ell_2$ as a
quotient~\cite{gasparis}, so they are not superprojective because the
corresponding quotient map is improjective (Proposition~\ref{superproj-char}).

\begin{Prop}\label{unconditional}
Every hereditarily-$\ell_1$ Banach space with a (countable or uncountable)
unconditional basis is subprojective.
\end{Prop}

\begin{Proof}
It was proved in \cite[Theorems 1 and~1a]{finol-wojtowicz} that every copy
of~$\ell_1$ in a Banach space~$X$ with a countable or uncountable
unconditional basis contains another copy of~$\ell_1$ which is complemented
in~$X$.
\end{Proof}

Later we will show a hereditarily-$\ell_1$ space $X_{bd}$ with a Schauder
basis which is not subprojective, so we cannot remove the unconditionality
condition in Proposition~\ref{unconditional}.

\medskip

We already know that
$C([0,\lambda],X)$ is subprojective when $X$ is subprojective and
$\lambda<\omega_1$ \cite{oikhberg-spinu}.
Next we improve this result.

\begin{Thm}
Let $X$ be a subprojective Banach space and let $\lambda$ be an arbitrary
ordinal. Then $C([0,\lambda],X)$ is subprojective.
\end{Thm}

\def\ColX{C_0([0,\lambda],X)}

\begin{Proof}
Observe that $C([0,\lambda],X) \equiv \ColX \oplus X$, so
$C([0,\lambda],X)$ is subprojective if and only if so is~$\ColX$
\cite[Proposition~2.2]{oikhberg-spinu}.
We will prove that $\ColX$ is subprojective by induction in~$\lambda$.

Assume that $C_0([0,\mu],X)$ is indeed subprojective for all $\mu<\lambda$.
If $\lambda$ is not a limit ordinal, then $\lambda = \mu+1$ for some~$\mu$
and $\ColX \equiv C_0([0,\mu],X)\oplus X$, and the result is clear.

Otherwise, if $\lambda$ is a limit ordinal, let $M$ be an
infinite-dimensional subspace of~$\ColX$ and define the projections
$$\map{P_\mu}{\ColX}{\ColX}$$
as $P_\mu(f) = f \chi_{[0,\mu]}$ for each $\mu<\lambda$.
If there exists $\mu<\lambda$ such that the restriction $P_\mu|_M$ is not
strictly singular, then there exists an infinite-dimensional subspace
$N\subseteq M$ such that $P_\mu|_N$ is an isomorphism. Since the range
of~$P_\mu$ is isometric to $C([0,\mu],X)$, which is subprojective by our
induction hypothesis, $N$ contains an infinite-dimensional subspace
complemented in $\ColX$ \cite[Corollary~2.4]{oikhberg-spinu}.

Assume now, on the other hand, that $P_\mu|_M$ is strictly singular for every
$\mu<\lambda$. We will construct a strictly increasing sequence of ordinals
$\lambda_1 < \lambda_2 < \ldots$ and a sequence $(f_n)_{n\in\N}$ of
normalised functions in~$M$ such that $\|P_{\lambda_{k-1}}(f_k)\| < 2^{-k}/8$
and $\|P_{\lambda_k}(f_k) - f_k\| < 2^{-k}/8$ for every $k\in\N$, where
we write $\lambda_0 = 0$ for convenience. To this end, let $k\in\N$,
and assume that $\lambda_{k-1}$ has already been obtained. By hypothesis,
$P_{\lambda_{k-1}}|_M$ is not an isomorphism, so there exists $f_k\in M$
such that $\|f_k\| = 1$ and $\|P_{\lambda_{k-1}}(f_k)\| < 2^{-k}/8$, and
then there is $\lambda_k \in (\lambda_{k-1},\lambda)$ such that
$\|P_{\lambda_k}(f_k) - f_k\| < 2^{-k}/8$, which finishes the inductive
construction process. Let $F = [f_k : k\in\N] \subseteq M$, which is
infinite-dimensional, and define the intervals $I_k =
(\lambda_{k-1},\lambda_k]$ and the operators $T_k = P_{\lambda_k} -
P_{\lambda_{k-1}}$ for every $k\in\N$; then $T_k(f) = f \chi_{I_k}$,
so each~$T_k$ is a norm-one projection and $T_iT_j = 0$ if $i\ne j$.

Let now $g_k = T_k(f_k) = P_{\lambda_k}(f_k) - P_{\lambda_{k-1}}(f_k)$
for each $k\in\N$; then
$$\|g_k - f_k\|
  \leq \|P_{\lambda_k}(f_k) - f_k\| + \|P_{\lambda_{k-1}}(f_k)\|
  < 2^{-k}/4,$$
so $1/2 < \|g_k\| < 3/2$ for every $k\in\N$. Note that
$C_0([0,\lambda])^* = \ell_1([0,\lambda))$ \cite[Theorem~14.24]{fabian-et-al}
and $C_0([0,\lambda],X)^* = (C_0([0,\lambda]) \itensor X)^* =
C_0([0,\lambda])^* \ptensor X^*$ \cite[Theorem~5.33]{ryan-book}, so
$$C_0([0,\lambda],X)^* = \ell_1([0,\lambda)) \ptensor X^* =
  \ell_1([0,\lambda),X^*)$$
and for each $k\in\N$ we can take $x_k \in \ColX^*$
with norm $\|x_k\| < 2$ such that $x_k(g_k) = 1$ and $x_k$ is
concentrated on~$I_k$, which makes $(g_n, x_n)_{n\in\N}$ a biorthogonal
sequence in $(\ColX,\ColX^*)$. In the spirit of the principle of small
perturbations \cite{bessaga-pelczynski}, let $K$ 
be the operator defined on $\ColX$ as
$K(f) = \sum_{n=1}^\infty x_n(f) (f_n - g_n)$; then
$$\sum_{n=1}^\infty \|x_n\| \|f_n - g_n\|
  < \sum_{n=1}^\infty 2^{-n}/2 = 1/2,$$
so $K$ is well defined and $U = I + K$ is an isomorphism on~$X$ that maps
$U(g_k) = f_k$ for every $k\in\N$.
Let $G = [g_k : k\in\N]$; then $U(G) = F$ and $G$ is infinite-dimensional.

We will now check that the supremum of the sequence $(\lambda_k)_{k\in\N}$
is $\lambda$ itself. Assume, to the contrary, that there existed some
$\mu < \lambda$ such that $\lambda_k\leq\mu$ for every $k\in\N$.
Then, for every $k\in\N$, we would have $P_\mu T_k =
P_\mu (P_{\lambda_k} - P_{\lambda_{k-1}}) =
(P_{\lambda_k} - P_{\lambda_{k-1}}) = T_k$, so $P_\mu(g_k) = g_k$ and
$P_\mu$ would be an isomorphism on~$G$. But then
$P_\mu U^{-1}$ would be an isomorphism on~$F$, where $U^{-1} = I - U^{-1} K$
is a compact perturbation of the identity, so $P_\mu$ would be upper
semi-Fredholm on~$F\subseteq M$, contradicting our assumption that
$P_\mu|_M$ is strictly singular.

This means, in turn, that $(x_n(f))_{n\in\N}$ is a null sequence
for every $f\in\ColX$, because each $x_k$ is supported on~$I_k$ and
$\|x_k\| < 2$, and we can define a projection
$$\map Q{\ColX}{\ColX}$$
as $Q(f)(\gamma) = x_k(f) g_k(\gamma)$ if $\gamma \in I_k$,
whose range is clearly~$G$.
Then $G$ is complemented in $\ColX$, and then so is $U^{-1}(G) = F
\subseteq M$, which proves that $\ColX$ is subprojective in this case too.
\end{Proof}

A Banach space $X$ has the Schur property when every weakly convergent
sequence in~$X$ is convergent.
Bourgain and Delbaen \cite{bourgain-delbaen} obtained two separable
$\mathcal{L}_\infty$-spaces $X_{bd}$ and~$Y_{bd}$ which admit Schauder bases
and satisfy the following properties:
\begin{itemize}
\item $X_{bd}$ has the Schur property, hence it is hereditarily~$\ell_1$; and
\item $Y_{bd}$ is hereditarily reflexive and $Y^*_{bd}$ is isomorphic
to~$\ell_1$.
\end{itemize}

To study these spaces, we need the following folklore result.
We include a proof for the convenience of the reader.

\begin{Prop}\label{c_0-quotient}
Every infinite-dimensional separable $\mathcal{L}_\infty$-space~$X$
has a quotient isomorphic to~$c_0$.
\end{Prop}

\begin{Proof}
Note that $X^*$ contains a sequence $(x^*_n)_{n\in\N}$ equivalent to the unit
vector basis of~$\ell_1$. Since $X$ is separable, passing to a subsequence we
can assume that $(x^*_n)_{n\in\N}$ is weak$^*$-convergent and, subtracting
the limit, that $(x^*_n)_{n\in\N}$ is weak$^*$-null.

We consider the operator $\map{T}{X}{c_0}$ defined as
$T(x) = (x^*_n(x))_{n\in\N}$. Since its conjugate operator~$T^*$ takes
the unit vector basis of~$\ell_1$ to the sequence $(x^*_n)_{n\in\N}$,
$T^*$ is an injection, hence $T$ is a surjection.
\end{Proof}

In Proposition~\ref{c_0-quotient} we can replace ``$X$ separable'' by
``the unit ball of $X^*$ is weak$^*$ sequentially compact''
\cite[Chapter XIII]{diestel}.

The next result for~$X_{bd}$ shows that an analogue of
Proposition~\ref{unconditional} for hereditarily-$\ell_1$ spaces is not valid
without further hypothesis.

\begin{Prop}
The spaces $X_{bd}$ and $Y_{bd}$ are neither subprojective nor
superprojective.
\end{Prop}

\begin{Proof}
The spaces $X_{bd}$ and $Y_{bd}$ are not subprojective because $\ell_1$
or a reflexive space cannot contain an infinite-dimensional
$\mathcal{L}_\infty$-space, and being an $\mathcal{L}_\infty$-space
is inherited by complemented subspaces.

For the other part, Proposition~\ref{ell-1} implies that $X_{bd}$ is not
superprojective, and for~$Y_{bd}$ (and also for~$X_{bd}$) we can apply
Proposition~\ref{c_0-quotient} to obtain a surjection $\map{T}{Y_{bd}}{c_0}$.
The kernel of~$T$ cannot be contained in any infinite-codimensional
complemented subspace~$M$, because $T$ would be an isomorphism on the
complement of~$M$ and $Y_{bd}$ does not contain copies of~$c_0$.
\end{Proof}

Note that $Y_{bd}^*\simeq \ell_1$ is subprojective,
but $X_{bd}^*\simeq C([0,1])^*$ is not.

\end{section}

\begin{section}{Sufficient conditions for superprojectivity}
\label{sect-cond}


An operator $\map{T}{X}{Y}$ is said to be \emph{unconditionally converging}
if there is no subspace~$M$ of~$X$ isomorphic to~$c_0$ such that the
restriction $T|_M$ is an isomorphism.
We denote the sets of unconditionally converging and weakly compact operators
from $X$ into~$Y$ by $U(X,Y)$ and $W(X,Y)$, respectively.

\begin{Def}
A Banach space~$X$ has property~(V) if $U(X,Y)\subseteq W(X,Y)$ for every
Banach space~$Y$; i.e. if every non-weakly compact operator $\map{T}{X}{Y}$
is an isomorphism on a subspace of~$X$ isomorphic to~$c_0$.
\end{Def}

It is well known that $C(K)$ spaces have property~(V), and it is not
difficult to see that property~(V) is inherited by quotients.
Property~(V) relates to superprojectivity because of the following result.

\begin{Thm}\label{abstract}
Let $X$ be a Banach space with property~(V) such that $X^*$ is
hereditarily~$\ell_1$. Then $X$ is superprojective.
\end{Thm}

\begin{Proof}
Let $M$ be an infinite-codimensional subspace of~$X$. Then $(X/M)^*$
contains a copy of~$\ell_1$, so $X/M$ admits an infinite-dimensional
separable quotient. Indeed, either $X/M$ has a quotient isomorphic to~$c_0$
or it contains a copy of~$\ell_1$ \cite{gonzalez-onieva}, in which case it
has a quotient isomorphic to~$\ell_2$. By passing to that further quotient,
we can assume that $X/M$ itself is separable.
However, $X^*$ is hereditarily~$\ell_1$, so $X/M$ is not reflexive,
and the quotient map $Q_M$ is not weakly compact.
By property~(V), there exists a subspace~$A$ of~$X$ isomorphic
to~$c_0$ such that $Q_M|_A$ is an isomorphism, where $Q_M(A)\simeq c_0$
is complemented because $X/M$ is separable.
Then $X/M = Q_M(A) \oplus B$, hence $X = A \oplus Q_M^{-1}(B)$ and
$M \subseteq Q_M^{-1}(B)$, so $X$ is superprojective.
\end{Proof}

\begin{Rem}\upshape
In the proof of Theorem~\ref{abstract} we need $X^*$ to be
hereditarily~$\ell_1$ to ensure the existence of separable quotients.
If this fact can be guaranteed for other reasons (e.g., $X$ separable)
we can replace ``$X^*$ hereditarily-$\ell_1$'' by the weaker condition
``$X$ does not admit infinite-dimensional reflexive quotients''.
\end{Rem}

Following Pietsch \cite[3.2.7]{pietsch}, we define $Sp(U^{-1}\circ K)$ as
the class of spaces $X$ satisfying that $U(X,Y)\subseteq K(X,Y)$ for every
Banach space~$Y$. This class admits a characterisation in terms of
property~(V) and the Schur property.
Let us first state an auxiliary result.

\begin{Prop}\label{dual-schur}
Let $X$ be a Banach space. Then the following are equivalent:
\begin{enumerate}
\itemn1 $X^*$ has the Schur property;
\itemn2 $X$ has the DPP and contains no copies of~$\ell_1$;
\itemn3 $W(X,Y)\subseteq K(X,Y)$ for every Banach space~$Y$.
\end{enumerate}
\end{Prop}

\begin{Proof}
For the equivalence between \rnbox1 and~\rnbox2, we refer to
\cite[Theorem~3]{diestel-DPP}.

For~\rnbox3, assume that $X^*$ has the Schur property, and take
$T\in W(X,Y)$; then $T^*\in W(Y^*,X^*) = K(Y^*,X^*)$, hence $T\in K(X,Y)$.
Conversely, if there exists a weakly null sequence $(x^*_n)_{n\in\N}$
in~$X^*$ that is not norm null, then the operator $\map{T}{X}{c_0}$ given by
$T(x) = (x^*_n(x))_{n\in\N}$ is weakly compact but not compact.
\end{Proof}

\begin{Prop}\label{abstract-char}
A Banach space~$X$ belongs to $Sp(U^{-1}\circ K)$ if and only if
it has property~(V) and its dual~$X^*$ has the Schur property.
\end{Prop}

\begin{Proof}
Property~(V) for~$X$ is equivalent to $U(X,Y)\subseteq W(X,Y)$ for every~$Y$,
and $X^*$ being Schur is equivalent to $W(X,Y)\subseteq K(X,Y)$ for every~$Y$
by Proposition~\ref{dual-schur}, which gives the desired result.
\end{Proof}

\begin{Cor}\label{abstract2}
Every Banach space in $Sp(U^{-1}\circ K)$ is superprojective.
\end{Cor}

\begin{Proof}
It is enough to observe that spaces with the Schur property
are hereditarily~$\ell_1$ and apply Proposition~\ref{abstract-char}
and Theorem~\ref{abstract}.
\end{Proof}

\begin{Cor}
A Banach space whose dual is isometric to~$\ell_1(\Gamma)$ belongs to
$Sp(U^{-1}\circ K)$, hence it is superprojective.
\end{Cor}

\begin{Proof}
The dual~$\ell_1(\Gamma)$ has the Schur property, and the space itself has
property~(V) by \cite[Corollary]{johnson-zippin}.
\end{Proof}

Note that, when $K$ is scattered, $C(K)^*$ is isometric to $\ell_1(K)$
\cite[Theorem 14.24]{fabian-et-al}, and that the space $Y_{bd}$ shows
that in the previous Corollary we cannot replace ``dual isometric''
by ``dual isomorphic''.
\medskip

The next results highlight the interest of the class $Sp(U^{-1}\circ K)$
by showing its stability under quotients, $c_0$-sums and projective tensor
products.

\begin{Prop}
The class $Sp(U^{-1}\circ K)$ is stable under passing to quotients.
\end{Prop}

\begin{Proof}
Suppose that $X$ belongs to $Sp(U^{-1}\circ K)$ and $\map{Q}{X}{Z}$
is a surjective operator.
Given $T\in U(Z,Y)$ we have $TQ\in U(X,Y)$. Then $TQ\in K(X,Y)$, hence
$T\in K(Z,Y)$.
\end{Proof}

\begin{Prop}\label{c_0(X_n)}
Given a sequence $(X_n)_{n\in\N}$ of spaces in $Sp(U^{-1}\circ K)$, the space
$c_0(X_n) = \{\, (x_n)_{n\in\N} : x_n\in X_n,\; (\|x_n\|)_{n\in\N} \in c_0 \,\}$
belongs to $Sp(U^{-1}\circ K)$.
\end{Prop}

\begin{Proof}
In the case $X_n=X$ for all~$n$, it was proved by Cembranos
\cite[Teorema~2]{cembranos-81} that $c_0(X_n)$ has property~(V) when
each $X_n$ does, and the proof is valid when the spaces $X_n$ are different.
Moreover $c_0(X_n)^*\equiv \ell_1(X_n^*)$ has the Schur property when each
$X_n^*$ does.
\end{Proof}

\begin{Thm}\label{tensor}
If the spaces $X$ and~$Y$ belong to $Sp(U^{-1}\circ K)$, then so does
$X\ptensor Y$.
\end{Thm}

\begin{Proof}
This is a consequence of two stability results for projective tensor
products. Ryan \cite[Corollary 3.4]{ryan} proved that if $X^*$ and~$Y^*$
have the Schur property then $(X\ptensor Y)^*$ also has the Schur property.
Moreover, if $X^*$ is Schur then $X$ contains no copies of~$\ell_1$ by
Proposition~\ref{dual-schur}, so any bounded sequence in~$X$ must contain a
weakly Cauchy subsequence. Since weakly Cauchy sequences in~$Y^*$, which is
Schur, must converge, this means that $L(X,Y^*) = K(X,Y^*)$, and it follows
from a result of Emmanuele and Hensgen \cite[Theorem 2]{emmanuele-hensgen}
that if $X$ and~$Y$ have property~(V) and $L(X,Y^*) = K(X,Y^*)$,
then $X\ptensor Y$ has property~(V).
\end{Proof}

\begin{Cor}
Let $X_1,\ldots,X_n$ be spaces belonging to $Sp(U^{-1}\circ K)$.
Then $X_1\ptensor\cdots\ptensor X_n$ is superprojective.
\end{Cor}

Note that $c_0\ptensor c_0$ is not an $\mathcal{L}_\infty$-space because
$(c_0\ptensor c_0)^{**}$ fails the DPP
\cite[Corollary 11]{gonzalez-gutierrez},
and it was proved in~\cite{galego-samuel} that $C(K) \ptensor C(L)$ is
subprojective when $K$ and~$L$ are countable compact.

\medskip

We do not know if $C(K,X)$ is superprojective when $K$ is scattered and $X$
is superprojective, but the following result gives a partial positive answer.
Recall that a Banach space $X$ has \emph{property~(u)} when for every
weakly Cauchy sequence $(x_n)_{n\in\N}$ in~$X$ there exists a weakly
unconditionally Cauchy series $\sum_{i=1}^\infty y_i$ so that
$\bigl( x_n - \sum_{i=1}^n y_i \bigr)_{n\in\N}$ is weakly null.
Banach spaces with an unconditional basis have property~(u)
\cite[Theorem~3]{pelczynski}.

\begin{Prop}\label{prop-u}
Suppose that $K$ is a scattered compact and $X$ is a Banach space with
property~(u) such that $X^*$ has the Schur property.
Then $C(K,X)$ belongs to $Sp(U^{-1}\circ K)$, and so it is superprojective.
\end{Prop}

\begin{Proof}
$X$ contains no copies of $\ell_1$ by Proposition~\ref{dual-schur}.
Since $K$ is scattered and $X$ has property~(u),
$C(K,X)$ has property~(V) \cite[Theorem~3]{cembranos-et-al}.
Moreover $C(K,X)^* \equiv \ell_1(K,X^*)$ has the Schur property, hence
$C(K,X)$ belongs to $Sp(U^{-1}\circ K)$ and it is superprojective by
Theorem~\ref{abstract}.
\end{Proof}

Pe\l czy\'nski proved \cite[Proposition~2]{pelczynski-2} that
a Banach space with property~(u) and containing no copies of~$\ell_1$ has
property~(V), so the condition on $X$ in Proposition~\ref{prop-u} implies
$X\in Sp(U^{-1}\circ K)$.

\subsection{The Hagler space}
\label{subsect-JH}

In \cite{hagler}, a Banach space~$\JH$ is constructed such that $\JH$ is
separable and hereditarily~$c_0$ and $\JH^*$ is nonseparable and has the
Schur property, hence it is hereditarily~$\ell_1$.
$\JH$ also has property~(S), which is defined as follows.

\begin{Def}
A Banach space~$X$ has property~(S) if every weakly null, non-norm null
sequence in~$X$ has a subsequence equivalent to the unit vector basis
of~$c_0$.
\end{Def}

Note that $\JH$ is subprojective by Proposition~\ref{hered-c0}.
Also, $\JH^*$ is not separable, so $\JH$ cannot admit an unconditional basis.

\begin{Prop}\label{Hagler}
The space~$\JH$ belongs to $Sp(U^{-1}\circ K)$, hence $\JH$ and
$\JH\ptensor\JH$ are superprojective.
Moreover $\JH^*$ is subprojective.
\end{Prop}

\begin{Proof}
Let us first see that $\JH$ belongs to $Sp(U^{-1}\circ K)$.
Let $\map{T}{\JH}{Y}$ be a non-compact operator, and let
$(y_n)_{n\in\N}$ be a bounded sequence in~$\JH$ such that $(T(y_n))_{n\in\N}$
has no convergent subsequence. Since $\JH$ contains no copies of~$\ell_1$ and
has property~(S), passing to subsequences and taking
$u_n := y_{2n} - y_{2n-1}$, we can assume that $(u_n)_{n\in\N}$ is equivalent
to the unit vector basis of~$c_0$ and $(T(u_n))_{n\in\N}$ is a seminormalised
basic sequence, and then,
since $\sum_{n\in\N} T(u_n)$ is weakly unconditionally Cauchy,
the sequence $(T(u_n))_{n\in\N}$ is equivalent to the unit vector basis
of~$c_0$ \cite[Corollary~V.7]{diestel}. Thus $T|_{[u_n]}$ is an isomorphism,
hence $\JH$ belongs to $Sp(U^{-1}\circ K)$
and Corollary~\ref{abstract2} implies that $\JH$ is superprojective.

To see that $\JH^*$ is subprojective, let $M$ be a subspace of~$\JH^*$.
As $\JH$ is separable and $\JH^*$ is Schur, we can find a sequence
$(x^*_n)_{n\in\N}$ in~$M$ equivalent to the unit vector basis of~$\ell_1$
which is weak$^*$-convergent to some $x_0^*\in X^*$. Let $y^*_n := x^*_n-x^*$;
by a remark of Johnson and Rosenthal \cite[Lemma~3.1.19]{gonzalez-martinez}
we can find a bounded sequence $(y_n)_{n\in\N}$ in~$X$ such that
$y^*_k(y_l) = \delta_{kl}$. Passing to a subsequence we can assume that
$(y_n)_{n\in\N}$ is weakly Cauchy, hence $(y_{2n} - y_{2n-1})_{n\in\N}$
is weakly null. We denote $z^*_n = y^*_{2n}$ and $z_n = y_{2n} - y_{2n-1}$.
Since $\JH$ has property~(S), we can assume that $(z_n)_{n\in\N}$
is equivalent to the unit vector basis of~$c_0$.

We consider the operators $\map{A}{X}{c_0}$ and $\map{B}{c_0}{X}$ defined
by $A(x) = (z^*_n(x))_{n\in\N}$ and $Be_n = z_n$. Then $P = BA$ is a
projection on~$X$ and $R(P^*)\subseteq M + \langle x_0\rangle$,
so $M$ contains a subspace complemented in~$X^*$.
\end{Proof}

The dual $\JH^*$ is not superprojective because it contains~$\ell_1$.

\begin{Prop}\label{C(K,JH)-super}
Let $K$ be a scattered compact.\ Then both\break
$C(K,\JH)\equiv C(K)\itensor\JH$ and $C(K)\ptensor\JH$ belong
to $Sp(U^{-1}\circ K)$, hence they are superprojective.
\end{Prop}

\begin{Proof}
It was proved by Knaust and Odell \cite[Theorem 2.1]{knaust-odell} that
property~(S) implies property~(u).
Since $\JH^*$ has the Schur property, Proposition~\ref{prop-u} implies
$C(K,\JH)\in Sp(U^{-1}\circ K)$.

The result for $C(K)\ptensor\JH$ follows from Theorem~\ref{tensor}.
\end{Proof}

\subsection{The Schreier space}


The Schreier space~$S$ is defined as the space of all scalar sequences
$x = (x_i)_{i\in\N}$ satisfying
$$\|x\|_S :=
  \sup \Bigl\{\, \sum_{i=1}^p |x_{n_i}| : p\leq n_1<\cdots n_p \,\Bigr\}
  < \infty.
$$

It satisfies the following properties:
\begin{itemize}
  \item[(a)] The unit vector basis is an unconditional basis for~$S$.
  \item[(b)] $S$ is a subspace of $C(\omega^\omega)$; as such,
  it is hereditarily~$c_0$.
  \item[(c)] $S$ fails the DPP \cite[Comments after Theorem~5]{diestel-DPP}.
  Hence $S^*$ is not Schur (Proposition~\ref{dual-schur}) and $S$ does not
  belong to $Sp(U^{-1}\circ K)$.
\end{itemize}

\begin{Prop}
The space~$S$ is subprojective and superprojective, and its dual~$S^*$ is
subprojective but not superprojective.
\end{Prop}

\begin{Proof}
$S$ is subprojective by Proposition~\ref{hered-c0}.
It is also separable, admits no infinite-dimensional reflexive
quotients \cite[Theorem~B and Corollary~1.10]{odell}, contains no copies
of~$\ell_1$, and satisfies property~(u) because it has an unconditional
basis. Thus $S$ has property~(V), and Theorem~\ref{abstract} implies that
$S$ is superprojective.

Its dual space~$S^*$ has an unconditional basis and, since $S$ admits no
infinite-dimensional reflexive quotient, $S^*$ contains no reflexive
subspace. Thus $S^*$ is hereditarily~$\ell_1$, hence it is subprojective
by Proposition~\ref{unconditional} and it is not superprojective by
Proposition~\ref{ell-1}.
\end{Proof}

Note that $S\notin Sp(U^{-1}\circ K)$ because $S^*$ is not Schur.
This is confirmed by the following result.

\begin{Prop}
The projective tensor products $S\ptensor S$ and $S\ptensor \ell_p$
($1<p<\infty$) are not superprojective.
\end{Prop}

\begin{Proof}
The dual space of $S\ptensor S$ can be identified with $L(S,S^*)$.
By \cite[Corollary~3.5]{gonzalez-pello} it is enough to show that that
there is a non-compact operator in $L(S,S^*)$.

Given $x = (x_i)_{i\in\N} \in S$, we denote the decreasing rearrangement
of~$(|x_i|)_{i\in\N}$ by $x^d = (x^d_i)_{i\in\N}$.
Note that, for each $n\in\N$, $x^d_{n}+\cdots+x^d_{2{n}-1}\leq \|x^d\|_S$,
so $x^d_{2n-1}\leq \|x^d\|_S/n$ and
$$\|x\|^2_2 = \|x^d\|^2_2 \leq 2({\textstyle\sum} 1/n^2)\|x^d\|^2_S
  \leq 2({\textstyle\sum} 1/n^2)\|x\|^2_S,$$
which means that $S\subseteq\ell_2$ and the natural inclusion
$\map{J}{S}{\ell_2}$ is a bounded operator, and then $\map{J^*J}{S}{S^*}$
is not compact.

The proof for $S\ptensor \ell_p$ is similar.
\end{Proof}

Observe that the previous argument does not apply to $S\ptensor c_0$.
We do not know if $S\ptensor c_0$ is superprojective.

\subsection{The predual of the Lorentz spaces~$d(w,1)$}

Given $p\geq 1$ and a non-increasing sequence of positive numbers
$w = (w_n)_{n\in\N}$, we consider the space $d(w,p)$ of all sequences of
scalars $x = (a_i)_{i\in\N}$ for which
$$\|x\| =
  \sup \biggl( \sum_{n=1}^\infty |a_{\pi(n)}|^p w_n \biggr)^{1/p}
  < \infty,$$
where the supremum is taken over all permutations $\pi$ of~$\N$.
Then $d(w,p)$ endowed with~$\|\cdot\|$ is a Banach space
\cite[Section~3a]{LT-1}.
To exclude trivial cases ($\ell_p$ or $\ell_\infty$) and normalise the
vectors we assume that $\lim_n w_n=0$, $\sum_n w_n=\infty$ and $w_1=1$.
In this case $d(w,p)$ is called a Lorentz sequence space
\cite[Definition~4.e.1]{LT-1}.

The unit vector basis $(e_n)_{n\in\N}$ is a symmetric basis for $d(w,1)$
and its biorthogonal sequence $(e^*_n)_{n\in\N}$ is a symmetric basis for
the predual $d(w,1)_*$ of $d(w,1)$. In particular, $d(w,1)_*$ contains
no copies of~$\ell_1$.

\begin{Prop}
The space $d(w,1)$ is subprojective and its predual $d(w,1)_*$ is
superprojective.
\end{Prop}

\begin{Proof}
The space $d(w,1)$ is hereditarily~$\ell_1$ \cite[Proposition~4.e.3]{LT-1},
hence subprojective by Proposition~\ref{unconditional}.

Since $d(w,1)_*$ has an unconditional basis, it satisfies property~(u)
\cite[Theorem~3]{pelczynski}, and $d(w,1)_*$ does not contain copies
of~$\ell_1$ because $d(w,1)$ is separable.
Then $d(w,1)_*$ has property~(V) and Theorem~\ref{abstract} implies
that it is superprojective.
\end{Proof}

\begin{Prop}
The space $d(w,1)$ fails the Schur property, so\break
$d(w,1)_*\notin Sp(U^{-1}\circ K)$.
\end{Prop}

\begin{Proof}
Note that $(e_n)_{n\in\N}$ is a symmetric basis in~$d(w,1)$ and
$$\lim_{n\to\infty} \frac{\|e_1+\cdots+e_n\|}{n} =
  \lim_{n\to\infty} \frac{w_1+\cdots+w_n}{n} = 0,$$
so $(e_n)_{n\in\N}$ is a normalised weakly null sequence.
\end{Proof}

\end{section}

\end{document}